\documentclass[twoside, 11pt]{article}
\usepackage{geometry}
\usepackage{times}
\usepackage[utf8]{inputenc}
\usepackage[vietnamese,english]{babel}
\usepackage[displaymath, tightpage]{preview}
\usepackage[all, cmtip]{xy}

\usepackage{amsmath, bm}
\usepackage{amsxtra}
\usepackage{amscd}
\usepackage{amsthm}
\usepackage{amsfonts}
\usepackage{amssymb}
\usepackage{bbm}
\usepackage{mathtools}
\usepackage[bitstream-charter]{mathdesign}
\usepackage[T1]{fontenc}

\usepackage{fancyhdr}
\usepackage{sectsty}
\sectionfont{\fontsize{13}{15}\selectfont}

\usepackage{tikz}
\usepackage{tikz-cd}
\usetikzlibrary{matrix,calc,arrows}
\usepackage[new]{old-arrows}

\usepackage{hyperref}
\usepackage[scaled=0.8]{helvet} 

\newtheoremstyle{theorem}
{3mm}
{1mm}
{\normalfont\itshape}
{}
{\normalfont\scshape}
{.}
{.5em}
{\thmname{#1}\thmnumber{ #2}\thmnote{ (#3)}}
\theoremstyle{theorem}
\newtheorem{theorem}{Theorem}[section]
\newtheorem{lemma}[theorem]{Lemma}
\newtheorem{proposition}[theorem]{Proposition}
\newtheorem{corollary}[theorem]{Corollary}

\newtheoremstyle{definition}
{3mm}
{1mm}
{\normalfont\normalfont}
{}
{\normalfont\scshape}
{.}
{.5em}
{\thmname{#1}\thmnumber{ #2}\thmnote{ (#3)}}
\theoremstyle{definition}
\newtheorem{remark}[theorem]{Remark}

\newtheorem{example}[theorem]{Example}
\newcommand\dqd{\foreignlanguage{vietnamese}{Đinh Quý Dương}}

\title{Rank-2 wobbly bundles \\
 from special divisors on spectral curves}
\author{\dqd}
\date{}

\begin{document}
\maketitle
	
\begin{abstract}
We study rank-2 wobbly bundles on a Riemann surface $C$ of genus $g\geq 2$,
i.e. semi-stable bundles admitting nonzero nilpotent Higgs fields,  
in terms of direct images of line bundles on 
smooth spectral curves $\tilde{C} \overset{\pi}{\rightarrow} C$.
We give a sufficient condition for a semi-stable bundle $E$ to be wobbly:
$E$ is a twist of $\pi_\ast \left(\mathcal{O}_{\tilde{C}}(\tilde{D}) \right)$ 
where the norm of $\tilde{D}$ is a summand of the divisor of a quadratic differential on $C$. 
We sketch the proof of the necessary condition statement, 
namely all wobbly bundles can be characterised as such,
and discuss how certain singularities of the wobbly locus arise 
from the Brill-Noether loci of spectral curves.
\end{abstract}

\section{Introduction}
Let $C$ be a Riemann surface of genus $g \geq 2$, 
$G$ a complex reductive group and $\check{G}$ its Langlands dual.
Laumon called a $G$-bundle $E$ on $C$ 
\textit{very stable} if it does not admit 
nonzero nilpotent Higgs fields $\phi \in H^0(C, \mathfrak{g}_E \otimes K_C )$ \cite{Lau88}.
He also showed that very stable bundles exist and are stable. 
Donagi-Pantev \cite{DP09} later called a bundle $E$ \textit{wobbly} 
if it is stable but not very stable. 
It is Drinfeld, according to Laumon, 
who conjectured that the wobbly locus is a divisor in the moduli space 
$\mathrm{Bun}_G$ of semistable bundles. 
In this article, we will work with vector bundles and 
say that a vector bundle $E$ is wobbly 
if it is \textit{semi-stable but not very stable}.

Drinfeld and Laumon noted early on the significance of the wobbly locus 
in the geometric Langlands correspondence.
This correspondence, which was recently proved in \cite{GLC-proof},
implies a very concrete fact: associated to a $\check{G}$-local system  
on $C$ is a local system on the complement  
of the wobbly locus in $\mathrm{Bun}_G$,
i.e. the locus of very stable bundles.
This way of thinking about the geometric Langlands correspondence 
in particular motivates a detailed study of the wobbly locus 
\cite{DP09, DP12, Pal-Pauly, PPe19, DP22, HH22, Pal22, DPS24, Pe24}

\paragraph{$\mathbb{C}^\ast$-action and wobbly Higgs bundles}
To put the main results of our paper in a larger context, 
it is helpful to recall the generalised notion of \textit{wobbly Higgs bundle}  
defined by Hausel-Hitchin \cite{HH22} for $G = GL_n$ or $SL_n$.
Let $\mathcal{M}_H$ be the moduli space 
of semi-stable Higgs bundles $(E, \phi)$ with gauge group $G$ on $C$. 
There is a natural $\mathbb{C}^\ast$-action on $\mathcal{M}_H$ that scales 
the Higgs fields, i.e. $a. (E,\phi) \coloneqq (E, a \phi) \in \mathcal{M}_H$ 
for $a \in \mathbb{C}^\ast$.
For any $(E, \phi) \in \mathcal{M}_H$, the limit $\underset{a \rightarrow 0}{\lim}
a. (E, \phi)$ exists and is a $\mathbb{C}^\ast$-invariant point in $\mathcal{M}_H$.
We say that 
$\mathcal{E}_0 = (E_0, \phi_0) \in \mathcal{M}_H^{\mathbb{C}^\ast}$
is \textit{very stable} if it is the only nilpotent Higgs bundle in 
the corresponding upward flow 
$$ W_{\mathcal{E}_0}^+ = 
\left\{ [E, \phi] \in \mathcal{M}_H \mid 
\underset{a \rightarrow 0}{\lim} a. [E, \phi] = \mathcal{E}_0
\right\}, $$
and is \textit{wobbly} otherwise.
If $E_0$ is stable, then $\mathcal{E}_0 = (E_0, 0)$  
and the upward flow $W_{\mathcal{E}_0}^+$ coincides with the cotangent 
fiber $T_{E_0}^\ast\mathrm{Bun}_G 
= \left\{(E_0, \phi) \mid \phi \in H^0(C, \mathfrak{g}_{E_0} \otimes K_C) \right\} \subset \mathcal{M}_H$: 
in this case we recover the notion of wobbly bundles.

Our aim is to study wobbly bundles and their locus from the point of view of 
direct images $\pi_\ast \mathcal{L}$ of line bundles $\mathcal{L}$ on
smooth spectral curves $\tilde{C} \overset{\pi}{\rightarrow} C$.
In this paper, we will study the cases where 
$G = GL_2$ and $G = SL_2$.
If the underlying bundle $E$ of
a rank-2 Higgs bundle $(E, \phi)$ with spectral curves $\tilde{C}$ 
is not semi-stable,
then $E = \pi_\ast \mathcal{L}$ where 
\begin{equation}\label{spectral-line-nss}
\mathcal{L} \simeq  \pi^\ast(L_E) \otimes \mathcal{O}_{\tilde{C}}(\tilde{D})
\end{equation}
with $L_E \hookrightarrow E$ the unique destabilizing subbundle of $E$
and $\tilde{D}$ an effective divisor of degree $\deg(\tilde{D}) < 2g-2$ \cite{Hit87a}.
The following is a result by Hausel-Hitchin 
applied to these cases.

\begin{theorem}\label{thm-HH} \cite{HH22}
Let $\tilde{C} \overset{\pi}{\rightarrow} C$ be a smooth rank-2 spectral curve 
and $\mathcal{L}$ a line bundle on $\tilde{C}$.
Assume that the underlying bundle $E$ 
of the Higgs bundle $(E, \phi) = \pi_\ast \mathcal{L}$ 
is not semi-stable, 
with (unique) destabilizing subbundle $L_E \hookrightarrow E$.
Then 
$[E_0, \phi_0] = \underset{a \rightarrow 0}{\lim} a. [E,\phi]$
is very stable if and only if 
the divisor $\tilde{D}$ in \eqref{spectral-line-nss} is reduced.
\end{theorem}
	
Our first main result is the analogue of Theorem \ref{thm-HH} 
for the case where $E$ is semi-stable.
Let us say that an effective divisor $D$ on $C$ is $Q$-special 
if it is a summand of the divisor of a quadratic differential,
i.e. $D < \mathrm{div}(q)$ for some $q \in H^0(C, K_C^2)$;
otherwise, we say that $D$ is $Q$-generic.
For brevity, in this article, we shall also say that an effective divisor 
$\tilde{D}$ on $\tilde{C}$ is $Q$-special if its norm 
$D = \mathrm{Nm}(\tilde{D})$ is $Q$-special.
It turns out that the condition of $\tilde{D}$ being $Q$-special  
is the proper analogue of the condition of being non-reduced 
in Theorem \ref{thm-HH}.

\begin{theorem}\label{main-thm-analogue}
Let $\tilde{C} \overset{\pi}{\rightarrow} C$ be a smooth rank-2 spectral curve 
and $\mathcal{L}$ a line bundle on $\tilde{C}$.
Assume that $E = \pi_\ast \mathcal{L}$ is semi-stable.
Then $E$ is wobbly if and only if 
\begin{equation*}
	\mathcal{L} \simeq \pi^\ast(L) \otimes \mathcal{O}_{\tilde{C}}(\tilde{D})
\end{equation*}
where $L$ is a line bundle on $C$ and $\tilde{D}$ is a $Q$-special divisor.
Furthermore, if $\tilde{D}$ has no summand of the form $\pi^\ast(P)$
for some effective divisor $P$ on $C$ then 
$E$ admits an embedding $L \overset{i}{\hookrightarrow} E$ 
and there is a canonical bijection between the corresponding space of nilpotent Higgs fields 
\begin{equation}
	H^0\left( C, K_C L^2 \det(E) \right) \simeq 
	\left\{ \phi \in H^0(C, \mathrm{End}(E)\otimes K_C) \text{ nilpotent},
	\ker(\phi) = i(L) \right\}
\end{equation}
and 
$$ \left\{ q\in H^0(C, K_C^2)
\mid D < \mathrm{div}(q) \right\}. $$
\end{theorem}

\paragraph{Components of the wobbly locus}
Consider now the case $G = SL_2$, 
i.e. we consider Higgs bundles whose underlying bundles $E$ 
have fixed determinant.
Given a fixed line bundle $\Lambda$ on $C$,
denote by $\mathrm{Bun}_{2, \Lambda}$ the moduli space of semi-stable rank-2 bundles $E$ with $\det(E) = \Lambda$.
Without loss of generality, we only need to consider the case where 
$\lambda = \deg(\Lambda)$ is equal to $0$ or $1$.
Consider the locus $\mathcal{W}_k$ which is 
the closure in $\mathrm{Bun}_{2, \Lambda}$ of 
\begin{equation}\label{def-Wk}
	\mathcal{W}_k^0 = 
	\{ E \in \mathrm{Bun}_{2, \Lambda} 
	\mid \exists \text{ sub-line bundle } L \text{ of } E, 
	h^0(C, K_C L^2 \Lambda^{-1}) > 0, \deg(K_C L^2 \Lambda^{-1})) = k  \}. 
\end{equation}
It is rather straightforward to show that such a locus $\mathcal{W}_k$ 
consists of wobbly bundles
\footnote{$E$ admits nonzero nilpotent Higgs fields if and only if $E$ has a sub-line bundle $L$ such that $h^0(C, K_C L^2 \Lambda^{-1}) > 0$ 
	(cf. Lemma \ref{lem-nilp-condition}). }.
Note that $k \equiv \lambda$ mod $2$ and, 
since $E$ is semi-stable and admits nonzero nilpotent Higgs fields, 
we have  
\begin{equation}\label{bound-Wk}
\lambda \leq k \leq 2g-2 - \lambda.
\end{equation} 
The following result by Pal and Pauly \cite{Pal-Pauly} settled the Drinfeld's conjecture 
for $G = SL_2$.

\begin{theorem} \label{thm-PP-Drinfeld} \cite{Pal-Pauly} 
The wobbly locus $\mathcal{W}$ in $\mathrm{Bun}_{2, \Lambda}$
is a divisor.    
Furthermore, $\mathcal{W}$ decomposes into 
\begin{equation}\label{wobbly-components}
\mathcal{W} = 
\begin{cases}
&\mathcal{W}_\lambda \cup \mathcal{W}_{\lambda + 2} \cup \dots \cup \mathcal{W}_g  
\qquad \quad \text{if } g \equiv \lambda \text{ mod } 2 \\
&\mathcal{W}_\lambda \cup \mathcal{W}_{\lambda + 2} \cup \dots \cup \mathcal{W}_{g-1}
\qquad \text{if } g \equiv \lambda - 1 \text{ mod } 2
\end{cases}.
\end{equation}
All loci $\mathcal{W}_k$ for $1 \leq k \leq g$ are irreducible divisors, 
except for $\mathcal{W}_0$ in the case $\lambda = 0$ is an union of $2^{2g}$
divisors each of which is equivalent to the Theta divisor of $\mathrm{Bun}_{2, \Lambda}$.
The locus $\mathcal{W}_k$ with $k > g$ in particular is of codimension $>1$ 
and is contained in either $\mathcal{W}_g$ (or $\mathcal{W}_{g-1}$)
if $g \equiv \lambda \text{ mod } 2$
(if $g \equiv \lambda - 1 \text{ mod } 2$).
\end{theorem}

Note that given a line bundle $\mathcal{L}$ on
a smooth spectral curve $\tilde{C} \overset{\pi}{\rightarrow} C$,
tensoring with a line bundle on $C$ changes neither the stability 
nor the wobbliness/very-stability of $\pi_\ast \mathcal{L}$.
Therefore, if $\tilde{D}$ is an effective divisor on $\tilde{C}$
and $\pi_\ast(\tilde{D})$ is a wobbly bundle, 
there are in total $2^{2g}$ choices of twisting by a line bundle 
to produce a wobbly bundle $E$ with $\det(E) = \Lambda$.

\begin{proposition}\label{prop-main}
	Let $\tilde{C} \overset{\pi}{\rightarrow} C$ be a smooth $SL_2$ spectral curve 
	and $\tilde{D}$ a $Q$-special divisor on $\tilde{C}$ of degree $d \in [2g-2+\lambda, 4g-4-\lambda]$.
	Assume that $\pi_\ast \mathcal{O}_{\tilde{C}}(\tilde{D})$ is semi-stable.
	Then given a square-root $L$ of $K_C \Lambda \otimes \mathcal{O}_C(-D)$,
	the rank-2 bundle
	\begin{equation*}
		L \otimes \pi_\ast \mathcal{O}_{\tilde{C}}(\tilde{D}).
	\end{equation*}
	is a wobbly bundle of determinant $\Lambda$ 
	and is contained in 
	$\mathcal{W}_{4g-4-d} \subset \mathrm{Bun}_{2, \Lambda}$.
	In particular, if $\tilde{D}$ does not contain any summand of the form
	$\pi^\ast(P)$ for some effective divisor $P$ on $C$, 
	then $E \in \mathcal{W}^0_{4g-4-d}$. 
\end{proposition}

Proposition \ref{prop-main} gives a sufficient condition for wobbliness.
We will sketch a proof of the necessary condition statement:
namely all wobbly bundles
$E \in \mathcal{W}_k$ for $\lambda \leq k \leq 2g-2-\lambda$
can be obtained as in Proposition \ref{prop-main}.
This is essentially a corollary of a theorem by Donagi-Pantev 
regarding resolving the direct image map $\pi_\ast$.

We end this paper by discussing how Brill-Noether loco on $\tilde{C}$ 
produce singularities of the wobbly locus via the direct image maps.

\section{Characterizing rank-2 nilpotent Higgs fields}
In this paper we will consider Higgs budles 
with gauge group $GL_2$ and $SL_2$. 
Denote by $\mathcal{M}_{2, \lambda}$ the moduli space 
of semi-stable $GL_2$-Higgs bundles $(E, \phi)$ with $\deg(E) = \lambda$, 
and $\mathcal{M}_{2, \Lambda}$ the moduli space of semi-stable 
$SL_2$-Higgs bundles $(E, \phi)$ where $\det(E) = \Lambda$ 
is a fixed line bundle on $C$.  
\subsection{Spectral line bundles in terms of divisors on spectral curves}
\paragraph{Spectral line bundles}
Let $\tilde{C} \overset{\pi}{\rightarrow} C$ be a smooth rank-2 spectral curve.
Locally, this means $\tilde{C} \subset T^\ast C$ is defined 
by $p_{\tilde{C}} = v^2 + b_1(z) v + b_2(z) = 0$, 
where $v$ is the cotangent fiber of $T^\ast C$ and
$(b_1, b_2)$ is an element of 
the Hitchin base $B = H^0(C, K_C) \oplus H^0(C, K_C^2)$.
Let $(E, \phi)$ be a Higgs bundle with 
$\tilde{C}$ as its associated spectral curve.
i.e. $p_{\tilde{C}}$ is the characteristic polynomial of $\phi$.
In case $\tilde{C}$ is smooth, its genus is $\tilde{g} = 4g-3$.

The cokernel of the map 
$$ \pi^\ast(E K_C^{-1}) \overset{\pi^\ast(\phi) - v}{\longrightarrow} \pi^\ast(E) $$
is a line bundle $\mathcal{L}_{(E,\phi)}$ on $\tilde{C}$ 
which we will call the spectral line bundle of $(E,\phi)$.
The spectral correspondence \cite{Hit87a, BNR} 
is a correspondence between Higgs bundles and their spectral line bundles,
namely that the underlying bundle can be recovered by the direct image 
$$ E \simeq \pi_\ast \mathcal{L}_{(E, \phi)},
\qquad \qquad \deg \left( \mathcal{L}_{(E,\phi)} \right) = 2g-2 +\lambda, $$
while the Higgs field $\phi$ is the direct image of 
$v$,
which coincides with the restriction of the canonical form on $T^\ast C$ 
restricted to $\tilde{C}$.
It follows that the generic fiber of the Hitchin fibration
$$ h: \mathcal{M}_{2, \lambda} \longrightarrow B = H^0(C, K_C) \oplus H^0(C, K_C^2) $$
is canonically isomorphic to $\mathrm{Pic}^{2g - 2 +\lambda}(\tilde{C})$  
which is a torsor over the Jacobian $J_{\tilde{C}}$ of $\tilde{C}$. 
For $G = SL_2$, 
the generic fiber of the Hitchin fibration 
$$h: \mathcal{M}_{2, \Lambda} \longrightarrow H^0(C, K_C^2) $$
is canonically isomorphic to
$$ \mathrm{Prym}^{2g - 2 +\lambda}(\tilde{C}, \Lambda) \coloneqq \{ \mathcal{L} \in \mathrm{Pic}^{2g - 2 +\lambda}(\tilde{C}) 
\mid \mathcal{L} \otimes \sigma^\ast\mathcal{L} \simeq \pi^\ast(\Lambda) \}  
\subset \mathrm{Pic}^{2g - 2 +\lambda}(\tilde{C}), $$
which is a torsor over the Prym variety 
$\mathrm{Prym}(\tilde{C}/C) = 
\{ \mathcal{L} \in J_{\tilde{C}} \mid \mathcal{L} \otimes \sigma^\ast \mathcal{L} 
\simeq \mathcal{O}_{\tilde{C}} \}  \subset J_{\tilde{C}}, $
where $\sigma$ is the involution of $\tilde{C}$.
The rational forgetful maps
$$ \mathrm{Pic}^{2g - 2 +\lambda}(\tilde{C}) \simeq h^{-1}(b) \dashrightarrow \mathrm{Bun}_{2, \lambda}, 
\qquad \qquad 
\mathrm{Prym}^{2g - 2 +\lambda}(\tilde{C}, \Lambda) \simeq h^{-1}(q) \dashrightarrow \mathrm{Bun}_{2, \Lambda}, $$ 
are given by taking direct images. 

\paragraph{Divisors induced by injection from line bundles}
Consider now an injection $i: L \rightarrow E$ from a line bundle $L$ into $E$,
which a priori is not necessarily an embedding and can vanish at some point on $C$.
Consider the composition 
\begin{equation} \label{def-BA-div}
	\pi^\ast(L) \overset{\pi^\ast i}{\longrightarrow} \pi^\ast(E) \longrightarrow \mathcal{L}_{(E, \phi)},
\end{equation}
which is a morphism between line bundles on $\tilde{C}$.
Denote by $\tilde{D}_i(\phi)$ the zero divisor of this morphism,
i.e. the support of its cokernel.
By definition, we have $\tilde{p} < \tilde{D}_i(\phi)$
if at $\tilde{p}$ the subspace $\pi^\ast(L)\mid_{\tilde{p}} < \pi^\ast(E)\mid_{\tilde{p}}$ 
is contained in the image of $\pi^\ast(\phi) - v$. 

It follows that we can write
$$\mathcal{L}_{(E, \phi)} \simeq \pi^\ast(L) \otimes \mathcal{O}_{\tilde{C}}(\tilde{D}_i(\phi)).$$
Note that if $L \overset{i}{\rightarrow} E$ vanishes at $p \in C$, then $\tilde{D}_i(\phi)$ 
contains a summand of the form $\pi^\ast(p)$, 
which is invariant under the involution $\sigma$ of $\tilde{C}$ 
that exchanges the two eigenvalues 
\footnote{If $p$ is a branch point of $C$, 
i.e. the two eigenvalues of $\phi\mid_p$ coincide, then $\pi^\ast(p) = 2 \pi^{-1}(p)$.}
of $\phi\mid_p$. 
On the base curve $C$,  consider the composition
\begin{equation}\label{map-on-C}
	c_i(\phi): L \overset{i}{\longrightarrow} E \overset{\phi}{\longrightarrow} E \otimes K_C 
	\longrightarrow L^{-1} \Lambda K_C,
	\qquad \qquad \Lambda = \det(E),
\end{equation}
where the last morphism is induced from the quotient of $i$.
We can regard $c_i(\phi)$ as a section of $K_C L^{-2} \Lambda$
and denote by $D_i(\phi)$ its zero divisor.

Denote by  
$$ \mathrm{Nm}: \mathrm{Pic}(\tilde{C}) \longrightarrow \mathrm{Pic}(C) $$
the norm map which maps a line bundle $\mathcal{O}_{\tilde{C}} \left( \sum_i n_i p_i \right)$ 
to the line bundle $\mathcal{O}_C \left( \sum_i n_i \pi(p_i) \right)$.
We will also use the notation $\mathrm{Nm}(\tilde{D})$ to denote the induced divisor on $C$ 
via the norm map.
\begin{proposition}\label{prop-Norm}
We have 
\begin{equation}
	D_i(\phi) = \mathrm{Nm}(\tilde{D}_i(\phi)).
\end{equation}
\end{proposition}
\begin{proof}
Suppose that $L \overset{i}{\hookrightarrow} E$ is an embedding 
and that in the local frames adapted to $i$ we have 
$\phi = 
\left( \begin{smallmatrix} a(z) & b(z) \\ c(z) & d(z) \end{smallmatrix} \right)$
where $c(z)$ is the local expression of $c_i(\phi)$.
Then at a point $\tilde{p} \in \tilde{C}$, the image of $\pi^\ast \phi - v$ 
which is a one-dimensional subspace of $\pi^\ast E \mid_{\tilde{p}}$,
coincides with $\pi^\ast L\mid_{\tilde{p}}$ 
if and only if $c(p) = 0$ for $p = \pi(\tilde{p})$ 
and $\tilde{p}$ is defined by the eigenvalue $v = d(z(p))$.
In other words, $\tilde{p} < \tilde{D}_i(\phi)$
if and only if $\pi(p) < D_i(\phi)$, 
and furthermore their multiplicities in the respective divisors are equal.
(Note that if $\tilde{p}$ is a ramification point of $\tilde{C} \rightarrow C$, 
the smoothness of $\tilde{C}$ requires $\tilde{p}$ and 
$p$ to both have multiplicity $1$ in $\tilde{D}_i(\phi)$ and $D_i(\phi)$ respectively.)
The case where $L \overset{i}{\rightarrow} E$ vanishes at a divisor $P$ 
follows by noting that $i$ is the composition of $L \overset{s_P}{\rightarrow} L(P)$ 
with an embedding $L(P) \hookrightarrow E$,
and that $\tilde{D}_i(\phi)$ in this case contains a summand of the form $\pi^\ast(P)$.
\end{proof}

\begin{proposition}\label{prop-inverse-BA}
Let $\tilde{D}$ be an effective divisor on a smooth spectral curve 
$\tilde{C} \overset{\pi}{\rightarrow} C$.
Let $L$ and $\Lambda$ be line bundles on $C$ satisfying 
$\mathcal{O}_C(D) \simeq K_C L^{-2} \Lambda$ 
where $D = \mathrm{Nm}(\tilde{D})$.
Then there exists a triple $(L \overset{i}{\rightarrow} E, \phi)$ 
unique up to scaling of $i$ such that $\tilde{D}_i(\phi) = \tilde{D}$
and $\det(E) = \Lambda$.
In particular, $i$ is an embedding if and only if $\tilde{D}$ contains 
no summand of the form $\pi^\ast(P)$ for some divisor $P$ on $C$.
\end{proposition}
\begin{proof}
Let $\mathcal{L} = \pi^\ast(L) \otimes \mathcal{O}_{\tilde{C}}(\tilde{D})$ 
and $(E, \phi) = \pi_\ast \mathcal{L}$.
By construction, $E$ admits an injection $i: L \rightarrow E$ 
defined by the morphism $\pi^\ast(L) \rightarrow \mathcal{L}$ 
that vanishes at $\tilde{D}$.
It follows that $\tilde{D} = \tilde{D}_i(\phi)$ 
by definition \eqref{def-BA-div}.
By Proposition \ref{prop-Norm} we can also calculate the determinant 
$\det(E) = K_C L^{-2} (-D)$.
Furthermore, by construction the injection $i$ 
vanishes at $P$ if and only if
$$\pi^\ast L \longrightarrow \pi^\ast E \longrightarrow \mathcal{L}_{(E, \phi)}$$
vanishes at $\pi^\ast(P)$,
i.e. $\pi^\ast(P) < \tilde{D}$.
\end{proof}

\begin{remark}
\begin{enumerate}
\item In \cite{DT23} these divisors $\tilde{D}_i(\phi)$ are called \textit{Baker-Akhiezer divisors}.
This is because, in local frames adapted to $\pi^\ast L$,
the line bundle $\pi^\ast(K_C^{-1}) \mathcal{L}_{(E,\phi)}$ 
of eigenvectors of $\pi^\ast(\phi)$ 
has sections of the form 
$\left(\begin{smallmatrix} v - d(z) \\ c(z)
\end{smallmatrix} \right)$,
reminiscent of the Baker-Akhiezer function in the theory of integrable systems.
\item Our definition \eqref{def-BA-div} of divisors $\tilde{D}_i(\phi)$ 
induced by triples $(L \overset{i}{\rightarrow} E, \phi)$ generalizes straightforwardly 
to higher rank cases.
In \cite{DDX}, we study the higher rank analogues of these divisors 
and the corresponding Lagrangian subspaces in the moduli spaces of Higgs bundles and 
flat connections. 
\end{enumerate}
\end{remark}

\subsection{Wobbly bundles as direct images of line bundles on spectral curves}
\paragraph{Kernels of nilpotent Higgs fields}
The following lemma was first noted in \cite{Pal17}
\begin{lemma}\label{lem-nilp-condition}
	Let $E$ be a rank-2 holomorphic vector bundle  on $C$ with $\det(E) = \Lambda$.
	Then $E$ admits nonzero nilpotent Higgs fields if and only if 
	it admits sub-line bundle $L \hookrightarrow E$ such that $h^0(C, K_C L^2 \Lambda^{-1}) > 0$.
\end{lemma}
\begin{proof}
	Suppose $E$ admits a nonzero nilpotent Higgs field $\phi_n$. Let $L = \ker(\phi_n)$ and consider the diagram
	\begin{equation*}
		\begin{tikzcd}
			&0 \arrow{r} &L \arrow{r} &E \arrow{r} \arrow{d}{\phi_n} &L^{-1} \Lambda \arrow{r} &0 \\
			&0 \arrow{r} &L\otimes K_C \arrow{r} &E \otimes K_C \arrow{r} \arrow{d}{\phi_n \otimes \mathbb{1}} &L^{-1} \Lambda K_C \arrow{r} &0 \\
			&  & &E \otimes K_C^2  &  &
		\end{tikzcd} \hspace{5pt} .
	\end{equation*}
	Since $L^{-1}\Lambda$ is isomorphic to the image of $\phi_n: E \rightarrow E\otimes K_C$, 
	there must be an embedding $L^{-1}\Lambda \hookrightarrow E \otimes K_C$. 
	But $(\phi_n \otimes \mathbb{1})\circ \phi_n = 0$, 
	so this embedding must factor through 
	$\ker(\phi_n \otimes \mathbb{1}) = L \otimes K_C \hookrightarrow E \otimes K_C$.
	This gives a non-trivial section of $K_C L^2 \Lambda^{-1}$.
	On the other hand, if $L$ is a sub-line bundle of $E$
	and $s$ is a non-trivial section of $K_C L^2 \Lambda^{-1}$,
	then the composition 
	$$ E \longrightarrow L^{-1}\Lambda \overset{s}{\longrightarrow} L \otimes K_C \longrightarrow E \otimes K_C $$
	clearly defines a nonzero nilpotent Higgs fields $\phi_n$ with $\ker(\phi_n) = L$.
\end{proof}

In the following, we will use the same notation $\phi_n$ 
for a nonzero nilpotent Higgs field $E \rightarrow E\otimes K_C$
and the corresponding map $L^{-1} \Lambda \rightarrow L K_C$.

\begin{proposition}\label{prop-nilp-quad}
	Let $E$ be a rank-2 bundle 
	which admits nonzero nilpotent Higgs fields with kernel 
	$i: L \hookrightarrow E$. 
	Let $\phi$ be a Higgs field on $E$ with smooth rank-2 spectral curve 
	$\tilde{C} \overset{\pi}{\rightarrow} C$ and $\tilde{D}_i(\phi)$ the associated divisor on $\tilde{C}$ defined by \eqref{def-BA-div}.
	Then $D = \mathrm{Nm}(\tilde{D}_i(\phi))$ is $Q$-special.
\end{proposition}
\begin{proof}
	Denote by $\Lambda$ the determinant of $E$.
	Consider the composition of the maps between line bundles on $C$ defined in \eqref{map-on-C}
	with the map in Lemma \ref{lem-nilp-condition} associated to kernels of nonzero nilpotent Higgs fields,
	\begin{equation*}
	\begin{tikzcd}
		&L \arrow{r} \arrow{rrrd} &E \arrow{r} &E\otimes K_C \arrow{r} & L^{-1} \Lambda K_C \arrow{d}{\phi_n} \\
		& & & &L  K_C^2
	\end{tikzcd}.
	\end{equation*}
	This composition is a quadratic differential $q$ with divisor 
	\begin{equation}\label{sum-div-1}
		\mathrm{div}(q) =  D + \mathrm{div}(\phi_n).
	\end{equation}
\end{proof}

\begin{proposition}\label{prop-quad-nilp}
	Let $\tilde{C} \overset{\pi}{\rightarrow} C$ be a smooth rank-2 spectral curve,
	and $\tilde{D}$ an effective divisor on $\tilde{C}$.
	If $D = \mathrm{Nm}(\tilde{D})$ 
	is a summand of the zero divisor of some quadratic differential,
	then the rank-2 bundle 
	$\pi_\ast \left( \mathcal{O}_{\tilde{C}}(\tilde{D}) \right)$
	has nonzero nilpotent Higgs fields.
\end{proposition}
\begin{proof}
	Let $(E, \phi) = \pi_\ast \left( \mathcal{O}_{\tilde{C}}(\tilde{D}) \right)$ 
	and $\Lambda = \det(E)$. 
	Suppose $\tilde{D}$ does not contain any summand which is invariant under the involution of $\tilde{C}$. 
	By Proposition \ref{prop-inverse-BA}, 
	$D$ is the divisor of a section $c \in H^0(C, K_C \Lambda)$ 
	defined by an embedding $\mathcal{O}_C \hookrightarrow E$.
	Let $q$ be a quadratic differential with $D < \mathrm{div}(q)$,
	and $\frac{q}{s} \in H^0(C, K_C \Lambda^{-1})$ with
	\begin{equation}\label{sum-div-2}
	\mathrm{div}(q) = D + \mathrm{div}\left( \frac{q}{s} \right).
	\end{equation}
	It follows from Lemma \ref{lem-nilp-condition} that $\frac{q}{s}$ defines nonzero 
	nilpotent Higgs fields on $E$ with kernel $\mathcal{O}_C$.
	
	Suppose now $\tilde{D} = \tilde{D}' + \pi^\ast( P)$ 
	where $\tilde{D}'$ does not contain any involution-invariant summand	
	and $P$ some effective divisor on $C$.
	Then there exists an embedding 
	$\mathcal{O}_C(P) \overset{i'}{\hookrightarrow} E$ 
	for which $i = i' \circ s_P$ and $\tilde{D}_{i'}(\phi) = \tilde{D}'$.
	Since $\mathrm{Nm}(\tilde{D}') < \mathrm{Nm}(\tilde{D})$ 
	is the summand of some quadratic differential,
	it follows from the argument above that $E$ admits nonzero nilpotent 
	Higgs fields with kernel $\mathcal{O}_C(P)$.
\end{proof}

\begin{proof}[Proof of Theorem \ref{main-thm-analogue}.]
Let $E$ be a bundle on $C$.
Given a line bundle $L$, the twist $L \otimes E$ is wobbly  
if and only if $E$ is wobbly.
Combining Propositions \ref{prop-nilp-quad} and \ref{prop-quad-nilp}
we can conclude that $\pi^\ast \mathcal{L}$ is wobbly if and only if 
$D$ is a $Q$-special divisor.
The bijection 
\begin{equation*}
	\left\{ \phi \in H^0(C, \mathrm{End}(E)\otimes K_C) \text{ nilpotent},
	\ker(\phi) = i(L) \right\}
	\longleftrightarrow
	\left\{ q\in H^0(C, K_C^2)
	\mid D < \mathrm{div}(q) \right\}
\end{equation*}
follows from equations \eqref{sum-div-1} and \eqref{sum-div-2}
\end{proof}

\section{Wobbly components for $G=SL_2$}
The Drinfeld's conjecture states that the wobbly locus 
shall be a divisor in the moduli space of semi-stable bundles.
For the $G = SL_2$ case, Pal and Pauly proved the conjecture and 
a decomposition of the wobbly divisor into irreducible divisors \eqref{wobbly-components} \cite{Pal-Pauly}.
Recently, Pal also announced a proof for $G = GL_n$ \cite{Pal22}.

Let us recall the wobbly components in the $G = SL_2$ case 
before proving Proposition \ref{prop-main},
which concerns with characterising wobbly components   
by taking direct image from the corresponding Picard components of spectral curves.
Consider the moduli space $\mathrm{Bun}_{2, \Lambda}$ 
of rank-2 semi-stable bundles with fixed determinant $\Lambda$, 
where $\lambda = \deg(\Lambda)$ is either $0$ or $1$.
It follows from Lemma \ref{lem-nilp-condition} that one 
can decompose the wobbly locus 
$\mathcal{W} \subset \mathrm{Bun}_{2, \Lambda}$ into 
\begin{equation*}
	\mathcal{W} = \bigcup_{\substack{\lambda \leq k \leq 2g-2 -\lambda, \\ k \equiv \lambda \text{ mod } 2}} \mathcal{W}_k^0, 
\end{equation*}
where  
\begin{equation}\label{def-Wk}
	\mathcal{W}_k^0 =
	\{ E \in \mathrm{Bun}_{2, \Lambda} 
	\mid \exists \text{ sub-line bundle } L \text{ of } E, 
	h^0(C, K_C L^2 \Lambda^{-1}) > 0, \deg(K_C L^2 \Lambda^{-1})) = k  \}.
\end{equation}
(cf. \eqref{bound-Wk} for the bounds of $k$).
Let
\begin{equation}\label{closure-wobbly-loci}
	\mathcal{W}_k = \overline{\mathcal{W}_k^0},
	\qquad \qquad k = \lambda, \dots, 2g-2-\lambda,
\end{equation}
be the closure in $\mathrm{Bun}_{2, \Lambda}$ of $\mathcal{W}_k^0$.

\begin{remark}
It is rather straightforward to show that 
the loci $\mathcal{W}_k$ are of pure codimension $1$
if $k \leq g$
and of codimension greater than $1$ for $k > g$.
The key technical step to prove Drinfeld's conjecture in this case
and the decomposition \eqref{wobbly-components},
therefore, is to show that the wobbly components of codimension great than $1$ 
lie in the closure of those of codimension $1$.
To this end, it was shown in \cite{Pal-Pauly} that 
the union 
$$ \bigcup_{\substack{ g < j \leq 2g - 2 - \lambda, \\ j \equiv \lambda \text{ mod } 2}} \mathcal{W}_j $$
is contained in the irreducible divisor 
$\mathcal{W}_g$ ($\mathcal{W}_{g-1}$) 
for $g = \lambda$ mod $2$ 
($g = \lambda - 1$ mod $2$).
\end{remark}

\paragraph{Proof of Proposition \ref{prop-main}}
	The fact that the determinant of $E = L \otimes \pi_\ast \mathcal{O}_{\tilde{C}}(\tilde{D})$ 
	is $\Lambda$ follows from Proposition \ref{prop-Norm},.
	It follows from Theorem \ref{main-thm-analogue} 
	and in particular equation \eqref{sum-div-2} that $E$ 
	is a wobbly bundle admitting an injection $L \rightarrow E$ 
	such that $\deg(K_C L^2 \Lambda^{-1}) = 4g-4-d$.
	In particular, if $\tilde{D}$ does not contain any summand of the form
	$\pi^\ast(P)$ then $L \rightarrow E$ is nowhere vanishing 
	and hence an embedding,
	and by definition $E \in \mathcal{W}^0_{4g-4-d}$.
	\hfill $\qed$
	
\begin{remark}
Recall that the injection $L \rightarrow E$ vanishes at $P$  
if and only if $\pi^\ast(P) < \tilde{D}$ and 
$E$ admits nilpotent Higgs fields with kernel $L(P) \hookrightarrow E$.
In this case, $E$ is contained in $\mathcal{W}^0_{4g-4 - d + 2\deg(P)} \cap \mathcal{W}_{4g-4-d}$.
\end{remark} 

\paragraph{Resolution of the direct image map}
We now sketch the proof that 
all wobbly bundles can be obtained as direct images of line bundles 
from a smooth spectral curve $\tilde{C} \overset{\pi}{\rightarrow} C$.  
We start by recalling the following result in \cite{DP09}, 
concerning how the wobbly locus in the $G = GL_2$ case 
can be described in terms of resolving the rational direct image map.
Let $\mathrm{Bun}_{2, \lambda}$ be the moduli space of semi-stable bundles 
of degree $\lambda$, which we can assume to be either $0$ or $1$.
Given a smooth rank-2 spectral curve 
$\tilde{C} \overset{\pi}{\rightarrow} C$,
consider the direct image map  
$\pi_\ast: \mathrm{Pic}^{2g-2 + \lambda}(\tilde{C}) \dashrightarrow \mathrm{Bun}_{2, \lambda}$.

\begin{theorem}\label{thm-DP} \cite{DP09}
	Let $\tilde{C} \overset{\pi}{\rightarrow} C$ be a smooth rank-2 spectral curve.
	\begin{enumerate}
		\item The rational map 
		$\pi_\ast: \mathrm{Pic}^{2g-2 + \lambda}(\tilde{C}) \dashrightarrow \mathrm{Bun}_{2, \lambda}$
		can be resolved to a morphism 
		$\Pi: \tilde{P} \rightarrow \mathrm{Bun}_{2, \lambda}$ 
		by a sequence of blow-ups with smooth centers.
		\begin{equation}
			\begin{tikzcd}
				& &\tilde{P} \arrow{dr}{\Pi} \arrow{dl}[swap]{f} & \\
				&\mathrm{Pic}^{2g-2 + \lambda}(\tilde{C}) \arrow[dashed]{rr}[swap]{\pi_\ast} & 
				&\mathrm{Bun}_{2, \lambda}
			\end{tikzcd}\hspace{5pt}.
		\end{equation}
		\item The image of the exceptional divisors $\bm{D} \subset \tilde{P}$
		along $\Pi$ coincides with the wobbly locus
		$\mathcal{W}_{2,\lambda} \subset \mathrm{Bun}_{2, \lambda}$.
		\item Given a translation invariant line bundle $\mathcal{H}$ on $\mathrm{Pic}^{2g-2 + \lambda}(\tilde{C})$,
		and a twist $\mathcal{H}'$ by exceptional divisors of the pull-back of $\mathcal{H}$ to $\tilde{P}$,
		the direct image $\Pi_\ast \mathcal{H}'$ is a quasi-parabolic Higgs sheaf 
		on $\mathrm{Bun}_{2, \lambda}$ with polar divisor equal to the wobbly locus.  
	\end{enumerate}
\end{theorem}

More generally, one can resolve the rational map 
$\mathcal{M}_H(C, GL_n) \dashrightarrow \mathrm{Bun}_{2, \lambda}$
by blowing-up the locus of stable Higgs bundles $(E, \phi)$ 
with unstable underlying bundles $E$.
The image in $\mathrm{Bun}_{2, \lambda}$ 
of the exceptional divisors via the resolution map 
is called the \textit{shaky} locus. 
It was conjectured by Donagi-Pantev and proved by Peón-Nieto \cite{Pe24} 
that the shaky and wobbly locus coincide.
Part 2 of Theorem \ref{thm-DP} states that, 
in the $G = GL_2$ case, the shaky locus in particular can be defined by blowing-up 
one single generic Hitchin fiber.

The following implication of Theorem \ref{thm-DP} was communicated to me by Donagi and Pantev.
\begin{corollary}\label{surjective-wobbly}
The rational maps $\pi_\ast: \mathrm{Pic}^{2g-2 + \lambda}(\tilde{C}) \dashrightarrow \mathrm{Bun}_{2, \lambda}$ is surjective.
\end{corollary}
\begin{proof}
That $\pi_\ast$ is surjective over the very stable locus follows 
from \cite{PPe19}.
Since $\mathrm{Bun}_{2, \lambda}$ is a Fano manifold of Picard number $1$,
the wobbly divisor $\mathcal{W} = \Pi(\bm{D})$ 
defines an ample line bundle $\mathcal{O}_{\mathrm{Bun}_{2,\lambda}}(\mathcal{W})$.
The pull-back along $\Pi$ of this line bundle is also ample.
The line bundle $\mathcal{O}_{\tilde{P}}(\bm{D})$, however, is not ample 
since $\bm{D}$ contracts along $f$.
It follows that $\Pi^{-1}(\mathcal{W})$ must contain $\bm{D}$ as a proper subset  
and that any point in $\mathcal{W} = \Pi(\bm{D})$ 
must be contained in the image of $\pi_\ast$.
\end{proof}

Theorem \ref{thm-DP} was announced for $G = GL_2$.
It is clear, however, that Parts $1$ and $2$ of the theorem 
are applicable also for $G = SL_2$.
This is because the direct image map $\pi_\ast$ is compatible with the passing 
between the two gauge groups:
given a smooth $SL_2$ spectral curve $\tilde{C}$, the diagram
\begin{equation}\label{gauge-groups-change}
	\begin{tikzcd}
		&\mathrm{Prym}^{2g-2 + \lambda}(\tilde{C}, \Lambda) \arrow[hook]{r} \arrow[dashed]{d}{\pi_\ast} 
		&\mathrm{Pic}^{2g-2 + \lambda}({\tilde{C}}) \arrow[dashed]{d}{\pi_\ast} 
		\arrow{r}{\simeq} &\left( \mathrm{Prym}^{2g-2 + \lambda}(\tilde{C}, \Lambda) \times J_C \right)/J_C[2] 
		\arrow[dashed]{d}{\pi_\ast} 
		\\
		&\mathrm{Bun}_{2, \Lambda} \arrow[hook]{r} &\mathrm{Bun}_{2, \lambda} 
		\arrow{r}{\simeq} &\left( \mathrm{Bun}_{2, \Lambda} \times J_C \right)/J_C[2]
	\end{tikzcd}\hspace{5pt}
\end{equation}
is commutative. 
Here $J_C[2]$ is the group of square-roots of the trivial line bundle 
$\mathcal{O}_C$, 
which acts naturally on the Cartesian products by component-wise tensoring 
(we identify $J_C$ with its pull-back $\pi^\ast J_C \subset J_{\tilde{C}})$.

\begin{theorem}\label{DP-SL2}
	Let $\tilde{C} \overset{\pi}{\rightarrow} C$ be an $SL_2$ spectral curve.
	\begin{enumerate}
		\item The rational map 
		$\pi_\ast: \mathrm{Prym}^{2g-2 + \lambda}(\tilde{C}, \Lambda) \dashrightarrow \mathrm{Bun}_{2, \Lambda}$ 
		can be resolved to a morphism $\Pi: \tilde{P}_2 \rightarrow \mathrm{Bun}_{2, \Lambda}$
		by a sequence of blowing up with smooth centers.
		\item The image of the exceptional divisors $\bm{D}_2 \subset \tilde{P}_2$
		along $\Pi$ coincides with the wobbly locus.
		\item The rational map 
		$\pi_\ast: \mathrm{Prym}^{2g-2 + \lambda}(\tilde{C}, \Lambda)  \dashrightarrow \mathrm{Bun}_{2, \Lambda}$ 
		is surjective.
	\end{enumerate}
\end{theorem}
\begin{proof}
	These are adaptations of Theorem \ref{thm-DP} and Corollary \ref{surjective-wobbly} to the case $G = SL_2$.
\end{proof}

\begin{remark}
It follows from Corollary \ref{surjective-wobbly} 
and Part 3 of Theorem \ref{DP-SL2} that, for $G = GL_2$ or $SL_2$,
the restriction $h_E$ of the Hitchin map 
to a cotangent fiber $T^\ast_E \mathrm{Bun}_G$ is at least surjective 
over the complement of the discrimant locus in the Hitchin base.
In particular, if $E$ is very stable then $h_E$ is surjective \cite{PPe19}.
We note that this result is relevant to the study of Brill-Noether loci of spectral curves
as discussed in \cite{N24}.
\end{remark}

\section{Singularities of the wobbly locus from Brill-Noether loci}
In the geometric Langlands correspondence,
associated to a local system on $C$
with gauge group being the Langlands dual $\check{G}$ 
is a local system on the complement of the wobbly locus 
in the moduli space $\mathrm{Bun}_G$ of semi-stable $G$-bundles on $C$.
In Donagi-Pantev's approach to the geometric Langlands correspondence,
these local systems conjecturally correspond via Mochizuki's non-abelian Hodge correspondence 
to parabolic Higgs bundles on $\mathrm{Bun}_G$,
which are direct images 
of line bundles on resolutions of the forgetful maps from Hitchin fibers to $\mathrm{Bun}_G$ \cite{DP09}.
Part 3 of Theorem \ref{thm-DP} in particular shows that, for $G = GL_2$, 
the singularities of these constructions are as expected.

To make use of Mochizuki's machinaries, 
one important criterion is to identify and resolve the singularities
of the wobbly locus 
that are more complicated than normal crossing and of codimension 1 in the wobbly locus $\mathcal{W}$
(codimension 2 in $\mathrm{Bun}_G$). 
A semi-stable bundle $E$ defines a smooth point of $\mathcal{W}$ if and only if 
it has exactly a one-dimensional space of nilpotent Higgs fields.
A typical normal-crossing singularity $E$ of $\mathcal{W}$ 
has two one-dimensional spaces $V_1, V_2 \subset \mathfrak{g}_E \otimes K_C$  
of nilpotent Higgs fields 
with $V_1 \cap V_2 = \{0\}$.
On the other hand, 
if $E$ has a space $V$ of nilpotent Higgs fields of dimension at least $2$, 
then $E$ is a singularity of the wobbly locus that is more complicated than  
normal-crossing.

\begin{remark}
	In the case $\lambda \equiv g -1 \text{ mod } 2$,
	the irreducible component $\mathcal{W}_{g+1}$ of the wobbly locus 
	is of codimension 1 in the wobbly locus (codimension 2 in $\mathrm{Bun}_{2, \Lambda}$).
	A bundle $E \in \mathcal{W}_{g+1}$ admits a maximal sub-line bundle $L$ 
	such that $\deg(K_C L^2 \Lambda^{-1}) = g+1$
	and hence there is a space $V \simeq H^0(C, K_C L^2 \Lambda^{-1})$ 
	of nilpotent Higgs fields on $E$ of dimension $\dim V \geq 2$.
	Hence $\mathcal{W}_{g+1}$ is an example of codimension-1 singularity in the wobbly 
	locus that is more complicated than normal-crossing.
\end{remark}

Given a smooth spectral curve $\tilde{C} \overset{\pi}{\rightarrow} C$, 
$r \geq 1$ and $d \geq 0$,
consider the corresponding Brill-Noether loci of $\tilde{C}$,
$$\mathrm{BN}^r_d(\tilde{C}) = 
\left\{ |\tilde{D}|  \in \mathrm{Pic}^d(\tilde{C}) \mid \dim|\tilde{D}|  \geq r-1  \right\}
= \left\{ \mathcal{O}_{\tilde{C}}(\tilde{D}) \in \mathrm{Pic}^d(\tilde{C}) \mid 
h^0\left( \tilde{C}, \mathcal{O}_{\tilde{C}}(\tilde{D}) \right) \geq r  \right\}, $$
where we have denoted by $|\tilde{D}|$ the linear series 
$\{ \tilde{D}' \text{ effective} \mid \tilde{D}' \sim \tilde{D} \}$.
The following corollary gives a lower bound of the dimension of nilpotent Higgs fields 
on a wobbly bundle $E \in \mathrm{Bun}_{2, \Lambda}$ 
obtained as in Proposition \ref{prop-main}
according to the Brill-Noether locus containing $|\tilde{D}|$.

\begin{proposition}\label{prop-BN}
Let $\tilde{C} \overset{\pi}{\rightarrow} C$ be a generic smooth $SL_2$ spectral curve, 
and $\tilde{D}$ be a $Q$-special divisor 
such that $|\tilde{D}| \in \mathrm{BN}^r_d(\tilde{C})$  
with $d = \lambda = \deg(\Lambda)$ mod $2$.
Suppose that $\pi_\ast\left( \mathcal{O}_{\tilde{C}}(\tilde{D}) \right)$
is semi-stable and that a generic divisor $\tilde{D}'$ 
in the linear series $|\tilde{D}|$ 
does not have any involution-invariant summand.
Then given a square-root $L$ of $K_C \Lambda \otimes \mathcal{O}_C(-D)$ 
where $D = \mathrm{Nm}(\tilde{D})$,
the wobbly bundle $E = L \otimes \pi_\ast \left( \mathcal{O}_{\tilde{C}} (\tilde{D}) \right)$
admits a space $V_L$ of nilpotent Higgs fields with kernels isomorphic to $L$
of dimension
\begin{equation}\label{lower-bound-nilp-Higgs}
	\dim V_L \geq 
	\begin{cases}
		r  &\hspace{20pt} \text{if }  3g-4+\lambda \leq d \leq 4g-4 - \lambda \\
		3g-4 -d + r &\hspace{20pt} \text{if }  d < 3g-4 
	\end{cases} \hspace{5pt}.
\end{equation}
\end{proposition}
\begin{proof}
Writing $\mathcal{L} = \pi^\ast(L) \otimes \mathcal{O}_{\tilde{C}} (\tilde{D})$,
by construction, we have a canonical isomorphism
$$\mathrm{Hom}(L, E) \overset{\simeq}{\longrightarrow} 
\mathrm{Hom}(\pi^\ast(L), \mathcal{L}) = 
H^0\left( \tilde{C}, \mathcal{O}_{\tilde{C}} (\tilde{D}) \right). $$
Hence $\dim \mathrm{Hom}(L, E) \geq r$ if $|\tilde{D}| \in \mathrm{BN}^r_d(\tilde{C})$.
By our hypothesis that a generic $\tilde{D}' \in |\tilde{D}|$ does not have an 
involution-invariant summand, 
there is an open dense subset $H' \subset \mathrm{Hom}(L, E)$ 
consisting of embeddings $L \hookrightarrow E$. 
The $\mathbb{C}^\ast$-scalings of such an embedding $i \in H'$
defines a subbundle in $E$ and a corresponding space of nilpotent Higgs fields 
$$V(i) = \left\{ \phi \in H^0(C, \mathrm{End}_0(E) \otimes K_C) \text{ nilpotent}, \ker(\phi) = i(L) \right\} \simeq H^0(C, K_C L^2 \Lambda^{-1}) $$
which has dimension
\begin{equation*}\label{lower-bound-nilp-Higgs-2}
	\dim V(i) \geq 
	\begin{cases}
		1  &\hspace{20pt} \text{if }  3g-4+\lambda \leq d \leq 4g-4 - \lambda \\
		3g-4 -d &\hspace{20pt} \text{if }  d < 3g-4 
	\end{cases} \hspace{5pt}.
\end{equation*}
The union
$$V_L = \bigcup_{i \in H'} V(i) 
\subset H^0(C, \mathrm{End}_0(E)\otimes K_C) $$
then is the space of nilpotent Higgs fields on $E$ with kernels 
isomorphic to $L$. Its dimension clearly has the lower bound given in \eqref{lower-bound-nilp-Higgs}.
\end{proof}

As we are interested in characterizing singularities of the wobbly locus,
let us note that if $|\tilde{D}| \in \mathrm{BN}^2_d$ 
then the wobbly bundle $E$ as constructed in Proposition \ref{prop-BN}
is a singularity of the wobbly locus that is more complicated than normal crossing.

\begin{corollary}\label{cor-sing-BN}
	Let $E$ be the wobbly bundle 
	$E = L\otimes \pi_\ast\left( \mathcal{O}_{\tilde{C}}(\tilde{D}) \right)$ as constructed in Proposition \ref{prop-BN} where  
	$|\tilde{D}| \in \mathrm{BN}^r_d(\tilde{C})$.
	Then $E$ is a singularity of the wobbly locus,
	which admits a space of nilpotent Higgs fields of dimension greater than $2$ 
	and is contained in the intersection 
	$\mathcal{W}_{4g-4-d} \cap \mathcal{W}_{4g-4-d'}$ for some $d' < d$.
\end{corollary}	
\begin{proof}
	It follows from \eqref{lower-bound-nilp-Higgs} that 
	$E$ admits a space of nilpotent Higgs fields of dimension greater than $1$ and hence 
	is a singularity of the wobbly locus. 
	Since $E$ is not split and 
	there are at least two embeddings $i_1, i_2: L \hookrightarrow E$ that are not 
	scalings of each other, 
	there exists a linear combination $i_3: L \rightarrow E$ of these two embeddings 
	such that $i_3$ vanishes at some effective divisor $P$.
	It follows from the argument in the proof of Proposition \ref{prop-main} that 
	$E$ admits nilpotent Higgs fields with kernel $L(P) \hookrightarrow E$, 
	and hence $E$ is contained in the intersection of wobbly components 
	$\mathcal{W}_{4g-4-d} \cap \mathcal{W}_{4g-4-d+2\deg(P)}.$
\end{proof}

We now would like to have an estimate of the dimension of 
the singular locus defined by such wobbly bundles with $r = 2$.
Let
$$ Q_d = \left\{ |\tilde{D} | \in \mathrm{Pic}^d(\tilde{C}) \mid 
\tilde{D} \text{ is } Q-\text{special} \right\} \subset \mathrm{BN}^1_d(\tilde{C}). $$
Note that for $d \leq 3g-4$ we have $Q_d = \mathrm{BN}^1_d(\tilde{C})$ and hence 
has dimension $d$, and for $3g-3 \leq d \leq 4g-4$ we have 
$\dim Q_d = 3g-4$.
On the other hand, the spectral curve $\tilde{C}$ 
a priori might not be sufficiently generic 
from the point of view of Brill-Noether theory, 
i.e. the dimension of $\mathrm{BN}^r_d(\tilde{C})$ is different from the naive estimate given by 
the Brill-Noether number 
$$\rho(\tilde{g}, r, d) = \tilde{g} - r(\tilde{g} - d +r -1), 
\qquad \qquad \tilde{g} = 4g-3, $$
in case $\rho(r, d)$ is non-negative. 
The actual locus of our interest, however, is 
the subset of $\mathrm{BN}^r_d(\tilde{C})$ 
whose direct images are stable bundles on $C$.
It is shown in \cite{N24} that for $g\geq 3$, the locus
$$\mathrm{BN}^2_{d, s}(\tilde{C}) = \left\{ |\tilde{D} | \in \mathrm{Pic}^2 (\tilde{C}) 
\mid h^0(\tilde{C}, \mathcal{O}_{\tilde{C}}(\tilde{D})  ) \geq 2, 
\pi_\ast  \mathcal{O}_{\tilde{C}}(\tilde{D})  \text{ stable}  \right\} 
\subset \mathrm{BN}^2_d(\tilde{C})$$ 
is empty for $d < 2g+1$ and has dimension $\rho(\tilde{g}, r, d)$ 
for $2g+1 \leq d \leq \tilde{g} = 4g-3$.
The singular locus of stable wobbly bundles as constructed in Corollary \ref{cor-sing-BN}
with $r=2$ is a twist of the direct images of the intersection 
$$ \bigcup_{2g + 2 -\lambda \leq d \leq 4g-4-\lambda} \left( \mathrm{BN}^2_{d,s}(\tilde{C}) \cap Q_d \right) \subset \mathrm{Pic}(\tilde{C}). $$
Denote this singular locus by $\mathcal{W}^2_{sing}$.

\begin{corollary}
Suppose that $g\geq 3$ and $Q_d$ 
is in general position in $\mathrm{BN}^1_d(\tilde{C})$.
Then $\mathcal{W}^2_{sing}$ has codimension 
$\mathrm{codim} \mathcal{W}^2_{sing} = 4$ 
($\mathrm{codim} \mathcal{W}^2_{sing} = 5$)
in $\mathrm{Bun}_{2, \Lambda}$ if $\lambda = 0$ ($\lambda = 1$ respectively).
\end{corollary}
\begin{proof}
Note that 
$$\rho(\tilde{g}, 2, d) = 2d - 4g +1. $$
For $3g-4 \leq d \leq 4g-4$,
$Q_d$ is a proper subset of $\mathrm{BN}^1_d(\tilde{C})$ of dimension $3g-4$,
and by our hypothesis is in general position.
Hence the intersection 
$Q_d \cap \mathrm{BN}^2_d(\tilde{C})$ has dimension 
$$\rho_{2, Q}(d) \coloneqq \rho(\tilde{g}, 2, d) + (3g-4) - d  = d - g - 3.$$ 
For the cases where $g = 3$ or $g = 4$, 
since $3g-4 = \dim \mathbb{P}H^0(C, K_C^2) < 2g+1$, 
we always need to take intersection of these two proper subsets 
in $\mathrm{BN}^1_d$. 
The dimension of $\mathcal{W}^2_{sing}$ is at most $\rho_{2, Q}(d'')$ where  
$$d'' = 
\begin{cases}
	4g-4 &\text{ if } \lambda = 0 \\
	4g-5 &\text{ if } \lambda = 1
\end{cases} \quad, $$
which gives codimension $4$ and $5$ in $\mathrm{Bun}_{2, \Lambda}$ respectively.
On the other hand, for $2g+1 \leq 3g-4$, i.e. $g \geq 5$, 
for any $2g+1 \leq d' \leq 3g-4$ we have 
$$ \rho(\tilde{g}, 2, d') \leq \rho(\tilde{g}, 2, 3g-4) = \rho_{2, Q}(3g-4) 
< \rho_{2,Q}(d'').$$
This again gives the codimension of $\mathcal{W}^2_{sing}$ as claimed.
\end{proof}

\begin{example}[Theta divisors]
In the moduli space $\mathrm{Bun}_{2, \mathcal{O}_C}$ of trivial determinant rank-2 bundles,
given a line bundle $L$ of degree $-g+1$,
one defines the associated generalised Theta divisor \cite{BNR}
$$ \Theta_{2, L} = 
\left\{ E \in \mathrm{Bun}_{2, \mathcal{O}_C} \mid h^0(C, L^{-1} E) \geq 1\right\}. $$
If we let $L$ be a square-root of $K_C^{-1}$, 
then $\Theta_{2, L}$ is an irreducible component in the wobbly divisor.
In fact, the wobbly component $\mathcal{W}_0$ is the union of all such divisors 
$$\mathcal{W}_0 = \bigcup_{L^2 \simeq K_C^{-1}} \Theta_{2, L}. $$
Given a smooth $SL_2$ spectral curve $\tilde{C}$,
consider its classical Theta divisor 
$$ \Theta_{\tilde{C}} 
\coloneqq \left\{ |\tilde{D}| \in \mathrm{Pic}^{\tilde{g} - 1}(\tilde{C}) 
\mid \tilde{D} \text{ effective} \right\} = \mathrm{BN}^1_{\tilde{g}-1}, 
\qquad \qquad \tilde{g} = 4g-3,$$
and an open dense subset $\Theta_{\tilde{C}, s} \subset \Theta_{\tilde{C}}$ consisting of line bundles 
whose direct images are stable.
Note that $\Theta_{\tilde{C}} = \mathrm{BN}^1_{\tilde{g}-1}$ 
contains $Q_{\tilde{g}-1}$.
Then the direct image of $Q_{\tilde{g}-1} \cap \Theta_{\tilde{C}, s}$ 
upon a twist by $L$ defines 
precisely the generalised Theta divisor $\Theta_{2, L}$.

The intersection of $Q_{\tilde{g}-1}$ with the 
higher Brill-Noether loci $\mathrm{BN}^r_{\tilde{g}-1}(\tilde{C})$, where $r \geq 2$,
$$ \mathrm{BN}^r_{\tilde{g}-1}(\tilde{C}) \cap Q_{\tilde{g}-1}$$ 
is non-empty and has dimension at least $3g- 3 - r^2$
provided that $r^2 \leq 3g-3$.
We see that points in these intersections,
which are the higher order vanishing loci of the Theta function on $J_{\tilde{C}}$, 
corresponds to singularities of the wobbly locus in Corollary \ref{cor-sing-BN}.
\end{example}

\vspace{10pt}	
\noindent \textbf{Acknowledgements:}
I would like to thank 
Ron Donagi, Rob Lazarsfeld, Federico Moretti, Tony Pantev, Vladimir Roubtsov and J\"org Teschner for helpful discussion.
Special thanks to Ron and Tony for explaining to me their result 
in \cite{DP09} and its implications.
This work is supported by funding and hospitality from 
Max Planck Institute for Mathematics (Bonn, Germany),
and by the German Research Foundation (DFG) under Germany’s Excellence Strategy – EXC 2121 “Quantum Universe” granted to the University of Hamburg.
The final revision of this paper was done with funding from 
DFG via the Walter Benjamin fellowship.

\vspace{10pt}
\noindent \textit{Contact:} duongdinh.mp@gmail.com 

\vspace{5pt}
\noindent 
\textit{Current addresses:} \\
Department of Mathtematics, University of Pennsylvania, Philadelphia, PA 19104, USA \\
University of Khánh Hoà, 1 Nguyen Chanh, Nha Trang, Khanh Hoa, Vietnam

\end{document}